
\magnification1440
\hfuzz5pt


\def\qed{\hfill{\vbox{\hrule\hbox{\vrule\kern3pt
                \vbox{\kern6pt}\kern3pt\vrule}\hrule}}}

\def\forc{{\,\parallel\joinrel\relbar\joinrel\relbar\,}}

\def\aa{{\alpha}}

\def\gg{{\gamma}}         
\def\kk{\kappa}

\def\oo{{\omega}}

\def\arrow{{\longrightarrow}}

\def\xx{\{x_0,x_1,x_2,x_3\}}

\smallskip
\centerline{\bf Two consistency results on set mappings }

\medskip
\centerline{\bf P\'eter Komj\'ath and Saharon Shelah} 

\bigskip
\noindent{\bf Abstract.} 
   It is consistent that there is a set mapping from 
   the four-tuples of $\omega_n$ into the finite subsets 
with no free subsets of size $t_n$ for some natural number $t_n$.
For any $n<\omega$ it is consistent that there is a set mapping from 
the pairs of $\omega_n$ into the finite subsets with no infinite 
free sets. 
For any $n<\omega$ it is consistent that there is a set mapping from 
the pairs of $\omega_n$ into $\omega_n$ with no uncountable 
free sets.
 
\bigskip
\noindent 
In\footnote{}{Research of the first author was partially supported by the 
              Hungarian National Science Research  Grant No.~T 016391. 
              Research of second the author was partially supported by 
              the United States-Israel Binational Science Foundation. 
              Publication number 645.}
this paper we consider some problems on 
{\sl set mappings}, that is, for our current purposes, 
functions of the type 
$f:[\kk]^k \to [\kk]^{<\mu}$  for some natural number $k$ and 
cardinals $\kk$, $\mu$, 
which satisfy 
$f(x)\cap x = \emptyset$ for $x\in [\kk]^k$. 
A subset $H$ of $\kk$ is called {\sl free} if 
$f(x)\cap H=\emptyset$ holds for every $x\in [H]^k$. 
The most central question of this area of combinatorial set theory is 
that given $k$, $\kk$, and $\mu$ how large free sets can be guaranteed. 
The investigation of the case $k=1$ was started in the thirties 
by Paul Tur\'an, who asked if there exists an infinite free 
set if $\mu=\oo$ and $\kk$ is the continuum. 
After G.~Gr\"unwald's affirmative answer ([4]) 
S.~Ruziewicz found the right conjecture ([10]); 
if $\kk>\mu$ then there is a free set of cardinal $\kk$ 
(remember, $k=1$ is assumed). 
Several cases were soon proved, 
for example S.~Piccard solved the case when $\kk$ is regular ([9]), 
but only in 1950 was the full conjecture established by Paul 
Erd\H os ([1]) with the assumption of GCH, and ten years later 
without this assumption, by A.~Hajnal ([5]). 
In the fifties Erd\H os and Hajnal started the research on the case 
$k>1$ following the observation of 
Kuratowski and Sierpi\'nski (see [4]) 
that for set mappings on $[\kappa]^k$ 
there always exists a free set of cardinal $k+1$ iff 
$\kk\geq \mu^{+k}$.

In ZFC alone,  Hajnal and M\'at\'e extended the 
Kuratowski-Sierpi\'nski results by showing ([6]) that if $k=2$ 
and $\kk\geq \mu^{+2}$ then there are arbitrarily large finite free 
sets, and Hajnal proved (see [3]) that a similar result holds 
for $k=3$, $\kk\geq \mu^{+3}$. 
One of the problems emphasized in [3] is if the result can be 
extended to $k=4$, $\kk\geq \mu^{+4}$. 
In Theorem 1 we show that it is not the case; 
for every natural number $n$ there exists a natural number 
$t_n$ such that for any given regular $\mu$ it is consistent 
that there is a set mapping 
$f:[\mu^{+n}]^4\to [\mu^{+n}]^{<\mu}$ with no free sets of size $t_n$. 
(We assume GCH in the ground model.)

As for the existence of infinite free sets, 
a special case of a theorem of Erd\H os and Hajnal states that 
under CH if $f:[\oo_2]^2\to [\oo_2]^{<\oo}$ is a set mapping then 
there is an uncountable free set for $f$ ([2]). 
Answering a question of [6] the first author proved that without CH even 
the existence of an infinite free set cannot be guaranteed [7]. 
Here we extend that result to arbitrary $\oo_n$. 
Using this result, we answer another question of Hajnal and M\'at\'e, 
by showing that it is consistent that  there exists a set mapping from 
the pairs of $\omega_n$ into  $[\omega_n]^1$ with no uncountable 
free sets. 

\medskip
Theorem 1 was proved by S.~Shelah; Theorem 2 and Corollary 3  
were subsequently proved by P.~Kom\-j\'ath. 

\medskip
\noindent{\bf Notation and Definitions.} 
We use the standard axiomatic set theory notation. 
Cardinals are identified with initial ordinals. 
If $S$ is a set and $\kk$ a cardinal, then 
$[S]^{\kk}= \{ X\subseteq S: |X| = \kk\}$, 
$[S]^{<\kk}= \{ X\subseteq S: |X| < \kk\}$, 
$[S]^{\leq \kk}= \{ X\subseteq S: |X| \leq \kk\}$. 
For $a$, $b$, $c$ and $r$ natural numbers, the Ramsey
symbol, $a\longrightarrow(b,c)^r$, means that the following
statement is true. 
Whenever the $r$-element subsets of an $a$ element set 
are colored with two colors, say 0 and 1, then either there exists 
a $b$-element subset with all its $r$-tuples colored 0 or there exists 
a $c$-element subset with all its $r$-tuples colored 1. 
The existence of an appropriate $a$ for any given $b$, $c$, $r$ 
is guaranteed by Ramsey's theorem [3]. 

\medskip
\noindent{\bf Acknowledgment.} 
The authors are grateful to the referee for several useful  
remarks which improved the exposition considerably.

\medskip
To formulate the next result, set $t_0=5$, $t_1=7$, in general, 
$t_{n+1}$ is the least number such that 
$t_{n+1} \arrow (t_n,7)^5$. 

\medskip

\noindent {\bf Theorem 1. } (GCH) 
         {\sl Assume that $n<\oo$, $\kk=\tau^{+n}$ for some 
              regular cardinal $\tau$. 
              Then it is consistent that GCH holds below $\tau$, 
              $2^\tau=\kk$ if $n>0$, and there is a set mapping 
              $f:[\kk]^4 \arrow [\kk]^{<\tau}$ 
              with no free subset of cardinal $t_n$. }

\medskip 
\noindent 
{\bf Proof.} 
By induction on $n$. 
Our set mapping will satisfy the additional condition that 
$f(\xx)\subseteq (x_1,x_2)$
(the ordinal interval) for all  $x_0<x_1<x_2<x_3$.

The case $n=0$ is obvious, 
since we can take $f(\xx)=(x_1,x_2)$. 

Assume that $V$ is a model of set theory 
satisfying the Theorem for $n$, and for $\tau^+$ in place of $\tau$. 
That is, for $\mu\leq\tau$, $2^\mu=\mu^+$ holds, 
and there is a set mapping 
$F:[\kk]^4\arrow [\kk]^{\leq\tau}$ satisfying 
$F(\xx)\subseteq (x_1,x_2)$ with no free subset of cardinality 
$t_n$. 
We are going to force with a notion of forcing $(P,\leq)$ 
in which the conditions will be some pairs of the form 
$(s,g)$ with 
$s\in [\kk]^{<\tau}$, $g:[s]^4\arrow [s]^{<\tau}$ satisfying 
$g(u)\subseteq F(u)$ for $u\in [s]^4$. 
Not all pairs as above will be in $P$ but if 
$(s,g)$, $(s',g')$ are in $P$ then $(s',g')$ will extend 
$(s,g)$ (in notation $(s',g')\leq (s,g)$) 
iff $s'\supseteq s$ and $g=g' | [s]^4$. 

To describe the condition for $(s,g)\in P$ 
we introduce two more definitions. 
If $U$ is a subset of $\kk$  then we call $U$  {\sl $F$-closed}, if 
$x_2\in F(x_0,x_1,x_3,x_4)$ holds whenever 
$x_0<x_1<x_2<x_3<x_4$ are in $U$. 
If $U$ is a subset of $s$  then we call $U$ {\sl $g$-free}, if 
$x_2\notin g(\{x_0,x_1,x_3,x_4\})$ holds 
for all $x_0<x_1<x_2<x_3<x_4$  in $U$. 
Now put $(s,g)$ into $P$ just in case there is no 7-element 
subset of $s$ which is $F$-closed and $g$-free.

Having defined the notion of forcing 
$(P,\leq)$, we are going to show some properties of it. 

\medskip 
\noindent{\bf Claim 1.} 
         {\sl $(P,\leq)$ is $<\tau$-closed.}

\medskip
\noindent 
{\bf Proof of Claim.} 
Immediate from the finite character of the definition. 
\qed

\medskip 
\noindent{\bf Claim 2.} 
         {\sl $(P,\leq)$ is $\tau^+$-c.c.}

\medskip
\noindent 
{\bf Proof of Claim.} 
Assume that $p_\xi=(s_\xi,g_\xi)\in P$ for $\xi<\tau^+$. 
Using the $\Delta$-system lemma we can assume that 
$s_\xi= a\cup b_\xi$ for some disjoint sets 
$\{ a\}\cup\{ b_\xi :\xi<\tau^+\}$. 
For $\xx\in [a]^4$, 
$g_\xi(\xx)$ is a subset of $F(\xx)$ of cardinal $<\tau$. 
As $\bigl|F(\xx)\bigr|\leq \tau$, and $|a|<\tau$, we can assume, 
by $\tau^{<\tau}=\tau$, that 
$g_\xi | [a]^4$ is the same for $\xi<\tau^+$. 
We show that any two $p_\xi$, $p_{\xi'}$ are compatible. 
Set $q=(a\cup b_\xi\cup b_{\xi'},g)$ 
where $g\supseteq g_\xi$, $g_{\xi'}$ and if 
$$
\bigl\{x_0,x_1,x_2,x_3\bigr\} 
                   \in 
[a\cup b_\xi\cup b_{\xi'}]^4 - 
[a\cup b_\xi]^4 - [a\cup b_{\xi'}]^4 
$$
then set 
$$
g(\xx)=(a\cup b_\xi\cup b_{\xi'})\cap F(\xx).
$$

We have to show that $q\in P$, that is, there is no 7-element 
$F$-closed, $g$-free subset of $s$. 
Assume that $B$ is such a set. 
As $p_\xi$, $p_{\xi'}$ are conditions, 
$B\not\subseteq a\cup b_\xi$, $B\not\subseteq a\cup b_{\xi'}$. 
There are, therefore, $\eta_0\in B\cap b_\xi$, 
$\eta_1\in B\cap b_{\xi'}$.

An easy calculation shows that no matter what position 
$\eta_0$, $\eta_1$ occupy in $B$, there is a five-tuple 
$y_0<y_1<y_2<y_3<y_4$ in $B$ such that $\eta_0$, $\eta_1\in 
\{y_0,y_1,y_3,y_4\}$. 
(This is the point where the choice of 7 plays role.) 
We get, therefore, that 
$g(\{y_0,y_1,y_3,y_4\})=F(\{y_0,y_1,y_3,y_4\})\cap s \ni y_2$ so 
$B$ cannot be $F$-closed and  $g$-free. 
\qed

\medskip
Let $G\subseteq P$ be a generic subset of $P$. 
Set $S=\bigcup \{s:(s,g)\in G\}$ and $f=\bigcup \{g:(s,g)\in G\}$. 
Clearly, $f$ is a set mapping of the required type on the set $S$. 

\medskip 
\noindent{\bf Claim 3.} 
         {\sl There is a $p\in P$ forcing that $|S|=\kk$.}

\medskip
\noindent 
{\bf Proof of Claim.} 
Otherwise, $1$ forces that $S$ is bounded in $\kk$, and 
as $(P,\leq)$ is $<\kk$-c.c., it forces a bound, say $\xi<\kk$. 
But as $\bigl(\{\xi\},\emptyset\bigr)\forc \xi\in S$, we get a 
contradiction. 
\qed

\medskip 
\noindent{\bf Claim 4.} 
         {\sl In $V[G]$, 
         $f$ has no free subset of cardinality $t_{n+1}$.}

\medskip
\noindent 
{\bf Proof of Claim.} 
Assume that $A\subseteq S$ is a  free subset of cardinality $t_{n+1}$. 
Color the five-tuples of $A$ as follows. 
If $\{x_0,x_1,x_2,x_3,x_4\}\in [A]^5$, $x_0<x_1<x_2<x_3<x_4$ and 
$x_2\in F(\{x_0,x_1,x_3,x_4\})$ then color 
$\{x_0,x_1,x_2,x_3,x_4\}$ 
by 1, otherwise by 0. 
As $t_{n+1} \arrow (t_n,7)^5$ either there is a homogeneous subset in color 
1 of cardinal 7 or there is a homogeneous subset of color 0 of size 
$t_n$. 
This latter possibility is excluded by the hypothesis on $F$ so 
we have the former. 
But that gives a 7-element subset which is $F$-closed and $f$-free 
and this is obviously excluded by the forcing. 
\qed

\medskip
Now Theorem 1 follows from the claims above by induction on $n$.
\qed

\bigskip
\noindent{\bf Theorem 2.} 
         (GCH) 
         {\sl If $\tau$ is a regular cardinal, $\kk<\tau^{+\oo}$, 
         then it is consistent that there is a set mapping 
         $f:[\kk]^2 \to [\kk]^{<\tau}$ with no infinite 
         free sets.}

\medskip
\noindent{\bf Proof.}          
For $\kappa\le\tau$, we can simply take $f(\{x,y\})=x$.

We are going to show, by induction on positive $n<\omega$
that it is  consistent that there exists for $\kappa=\tau^{+n}$
a set mapping $f$ on $[\kappa]^2$ as required.  
It will also satisfy $f(\{x,y\})\subseteq x$ for $x<y<\kappa$.

The case $n=1$ can also be proved in ZFC. 
If $x<\kk=\tau^+$, enumerate $x$ as $x=\{\gg_x(i):i<\tau\}$. 
If $x<y$ then let $i(x,y)$ be that index $i$ for which $x=\gg_y(i)$ 
holds. 
Now set 
$f\bigl(\{x,y\}\bigr)=\bigl\{\gg_x(i):i\leq i(x,y)\bigr\}$. 
If $x_0<x_1<\cdots$ are the elements of an infinite free set then 
$i(x_0,x_1) > i(x_1,x_2)> \cdots$ which is impossible. 

Assume now that $\tau<\kk$, GCH holds up to and including $\tau$ 
and there is a set mapping 
$F:[\kk]^2 \to [\kk]^{\leq\tau}$ with no infinite 
free sets and with $F\bigl(\{x,y\}\bigr)\subseteq x$ for $x<y<\kk$. 
We are going to define a $<\tau$-closed partial ordering 
$(P,\leq)$ which adds a set mapping 
$f:[S]^2 \to [S]^{<\tau}$ for some $S\in [\kk]^\kk$ 
and  with no infinite free sets. 
It will also satisfy 
$f(\{x,y\})\subseteq F(\{x,y\})$ for all $\{x,y\}\subseteq S$. 

An element of $P$ will be a triplet of the form 
$p=(s,g,r)$ where $s\in [\kk]^{<\tau}$, 
$g:[s]^2\to [s]^{<\tau}$ is a set mapping with $g\subseteq F$. 
If $U$ is a subset of $\kk$  then we call $U$  {\sl $F$-closed}, if 
$x\in F(y,z)$ holds if $x<y<z$ are in $U$. 
If $U$ is a subset of $s$  then we call $U$ {\sl $g$-free}, if 
$x\notin g(\{y,z\})$ holds for $x<y<z$  in $U$. 
We require that there be no infinite $g$-free, $F$-closed subsets  
of $s$ and $r$ will be a rank function witnessing this. 
For this, we call a finite subset $u\in [s]^{<\oo}$ {\sl secured} if 
$|u|\geq 3$, $u$ is $g$-free and $F$-closed. 
What we assume on $r$ is that it is a function from the secured 
subsets to $\tau$ with 
$r(u)>r(v)$ if $v$ properly end-extends $u$. 
$p'=(s',g',r')$ extends $p=(s,g,r)$ if $s' \subseteq s$, 
$g'\subseteq g$, $r'\subseteq r$. 

It is obvious that $(P,\leq)$ is transitive and $<\tau$-closed. 

\medskip
\noindent{\bf Claim 1.} 
         {\sl $(P,\leq)$ is $\tau^+$-c.c.} 

\medskip
\noindent 
{\bf Proof of Claim.} 
Assume, for a contradiction, that we are given 
$\tau^+$ conditions, 
$p_\xi=(s_\xi,g_\xi,r_\xi)\in P$ for $\xi <\tau^+$.          
By the $\Delta$-system lemma we can assume that there are disjoint sets 
$\{a\}\cup\{b_\xi:\xi<\tau^+\}$ such that $s_\xi=a\cup b_\xi$. 
As for $x$, $y\in a$, since 
$g_\xi(\{x,y\})\in [F(\{x,y\})]^{<\tau}$, by removing at most $\tau$ 
members from the family we can assume that 
$F\bigl(\{x,y\}\bigr)\cap b_\xi = \emptyset$ holds for $x$, $y\in a$. 
Then, $g_\xi\bigl(\{x,y\}\bigr)\subseteq a$, and with one more 
shrinking, we can assume that $g_\xi\bigl(\{x,y\}\bigr)$ is independent 
of $\xi$. 
We can also assume that the functions $r_\xi$ are identical on the 
secured subsets of $a$. 

Assume now that $\xi<\xi'<\tau^+$, we want to find a common extension of 
$p_\xi$ and $p_{\xi'}$. 
Set $q=(a\cup b_\xi \cup b_{\xi'},g,r)$ where $g\supseteq g_\xi\cup g_{\xi'}$ 
is the maximal extension, that is, 
$g\bigl(\{x,y\}\bigr) = (a\cup b_\xi \cup b_{\xi'})\cap 
F\bigl(\{x,y\}\bigr)$ if  
$\{x,y\}\cap b_\xi\ne\emptyset$ and
$\{x,y\}\cap b_{\xi'}\ne\emptyset$. 

We now consider if we can define $r$.  
As $q$ is the union of two conditions both omitting infinite 
$g$-free, $F$-closed sets, $q$ won't have such sets, either. 
So {\sl some}  rank function $r$ can be defined; the question is, 
if one extending  $r_\xi$, $r_{\xi'}$ can be given. 
To show this, it suffices to prove, that if $u$ is a $g$-free, 
$F$-closed set, which is new, that is, has points in 
$b_\xi$, as well as in $b_{\xi'}$, then it cannot end extend an ``old'' 
secured set (one in $p_\xi$ or in $p_{\xi'}$). 
Assume that $x_0<x_1<x_2<\cdots$ are the elements of $u$. 
If $x_i\in b_\xi$, $x_j\in b_{\xi'}$, and $i$, $j\neq 0$, then 
$x_0\in F\bigl(\{x_i,x_j\}\bigr)$, so 
$x_0\in g\bigl(\{x_i,x_j\}\bigr)$ by the definition of $g$ 
and so our set is not $g$-free. 
We get, therefore, that $x_0$ is the {\sl only} element of $u\cap b_\xi$ 
(say). 
The possibility that both $x_1$ and $x_2$ are in $a$ is ruled out by our 
above condition that $b_\xi\cap F\bigl(\{x_1,x_2\}\bigr)=\emptyset$. 
This means that $\{x_0,x_1,x_2\}$ is a ``new'' set, so $u$ is indeed not 
an end extension of an old secured set as we assumed that secured sets 
have at least three elements. 
\qed

\medskip
If $G\subseteq P$ is a generic subset, 
then define 
$S=\bigcup \{s:(s,g,r)\in G\}$, 
$f=\bigcup \{g:(s,g,r)\in G\}$, 
$R=\bigcup \{r:(s,g,r)\in G\}$, 
. 

\medskip
\noindent{\bf Claim 2.} {\sl $|S|=\kk$.} 

\medskip
\noindent 
{\bf Proof of Claim.} 
As in the corresponding proof in Theorem 1. 
\qed

\medskip
\noindent{\bf Claim 3.} {\sl $F$ has no infinite free set in $V[G]$.} 

\medskip
\noindent 
{\bf Proof of Claim.} 
This is a well-known fact. 
It follows from the rank characterization of the nonexistence 
of free sets. 
\qed

\medskip
\noindent{\bf Claim 4.} {\sl $f$ has no infinite free set.} 

\medskip
\noindent 
{\bf Proof of Claim.} 
Assume that $x_0<x_1<\cdots$ form an infinite $f$-free set. 
By Ramsey's theorem we can assume that either for every triplet  $i<j<k<\oo$, 
$x_i\in F(x_j,x_k)$ holds or for every triplet  $i<j<k<\oo$, 
$x_i\notin F(x_j,x_k)$ holds. 
The latter  is impossible by Claim 3. 
Therefore $\{x_0,x_1,\dots\}$ is $f$-free, $F$-closed, but then 
$R\bigl(\{x_0,x_1,x_2\}\bigr) > R\bigl(\{x_0,x_1,x_2,x_3\}\bigr) >\cdots$, 
which is impossible. 
\qed

\medskip
An easy application of Theorem 2 solves another problem of [6]. 

\medskip
\noindent{\bf Corollary 3.} 
         {\sl  For every $n<\oo$ it is consistent that 
               there exists a set mapping $f:[\oo_n]^2\to [\oo_n]^1$ 
               with no uncountable free set.}

\medskip
\noindent 
{\bf Proof.} 
Applying Theorem 2 assume that $F:[\oo_n]^2\to [\oo_n]^{\aleph_0}$ 
is a set mapping with no infinite free sets 
so that $F(\{x,y\})\subseteq x$
for all $x<y<\omega_n$. 
Define the notion of forcing as follows,  
$(s,g)\in P$ iff 
$s\in [\oo_n]^{<\oo}$, $g:[s]^2\to [s]^1$, 
and
$g(u)\subseteq F(u)$ for all $u\in[s]^2$.
Set $(s',g')\leq (s,g)$ iff $s'\supseteq s$, $g'\supseteq g$. 

\medskip
\noindent{\bf Claim 1.} {\sl If $\aa<\oo_n$ then the set 
                        $D_\aa=\{(s,g):\aa\in s\}$ is dense in 
                        $(P,\leq)$.} 

\medskip
\noindent 
{\bf Proof of Claim.} 
Straightforward. 
\qed

\medskip
\noindent{\bf Claim 2.} {\sl $(P,\leq)$ is c.c.c.} 

\medskip
\noindent 
{\bf Proof of Claim.} 
Assume that $p_\xi\in P$ for $\xi<\oo_1$. 
By the usual thinning out procedure we can assume that 
$p_\xi=(s\cup s_\xi,g_\xi)$  where $s_\xi\cap F(x,y)=\emptyset$ holds 
for $x$, $y\in s$, and the functions $g_\xi | [s]^2$ 
are identical. 
Now any two $p_{\xi}$-s are compatible. 
\qed

\medskip
If $G\subseteq P$ is a generic set, put 
$f=\bigcup\{g:(s,g)\in G\}$. 
\medskip
\noindent{\bf Claim 3.} {\sl $f$ has no uncountable free set.} 

\medskip
\noindent 
{\bf Proof of Claim.} 
Assume that $p\forc X$ is an uncountable free set. 
There are, for $\xi<\oo_1$, conditions $p_\xi\leq p$ 
and ordinals $\aa_{\xi}$ with $p_\xi\forc \aa_\xi \in X$. 
Again, we can assume, that 
$p_\xi= (s\cup s_\xi, g_\xi)$, $\aa_\xi \in s_\xi$, and the functions 
$g_\xi \cap [s]^2$ are identical. 
As $F$ has no infinite free sets (``no uncountable'' suffices) 
there are ordinals $\xi_0$, $\xi_1$, $\xi_2<\oo_1$ such that 
$\aa_{\xi_0}\in F\bigl(\aa_{\xi_1},\aa_{\xi_2} \bigr)$. 
We can now extend $p$ to a condition $p'=(s',g')$ where 
$$
s' = s \cup s_{\xi_0} \cup s_{\xi_1} \cup s_{\xi_2},
$$
$g'$ extends $g_{\xi_0}$, $g_{\xi_1}$, $g_{\xi_2}$ and 
$g'\bigl(\{\aa_{\xi_1},\aa_{\xi_2}\} \bigr) = \aa_{\xi_0}$. 
\qed

\medskip
Now Corollary 3 follows from the claims above.
\qed
         
\bigskip
\centerline{\bf References}

\medskip
\item{[1]} P.~Erd\H os: 
         Some remarks on set theory, 
         {\sl Proceedings of the American Mathematical Society}, 
         {\bf 1}(1950), 127--141. 

\item{[2]} P.~Erd\H os, A.~Hajnal: 
         On the structure of set mappings, 
         {\sl Acta Mathematica Acad.~Sci.Hung.}, 
         {\bf 9}(1958), 111--131.

\item{[3]} P.~Erd\H os, A.~Hajnal, A.~M\'at\'e, R.~Rado: 
          {\sl Combinatorial Set Theory: 
          Partition Relations for Cardinals}, 
          North-Holland, Akad\'emiai Kiad\'o, 1984. 

\item{[4]} G.~Gr\"unwald: 
         Egy halmazelm\'eleti t\'etelr\H ol, 
         {\sl Mathematikai \'es Fizikai Lapok}, 
         {\bf 44}(1937), 51--53. 

\item{[5]} A.~Hajnal: 
         Proof of a conjecture of S. Ruziewicz, {\sl Fund. Math. 
         } {\bf 50} (1961), 123-128.

\item{[6]} A.~Hajnal, A.~M\'at\'e: 
          Set mappings, partitions, and chromatic numbers, 
          in: {\sl Logic Colloquium '73}, Bristol, 
          North-Holland, 1975, 347--379. 
          
\item{[7]} P.~Komj\'ath: 
          A set mapping with no infinite free subsets, 
          {\sl Journal of Symbolic Logic \bf 56} (1991), 
          304--306.

\item{[8]} K.~Kuratowski: 
          Sur une charact\'erization des alephs, 
          {\sl Fundamen\-ta Ma\-thematicae}, {\bf 38}(1951), 
          14--17. 

\item{[9]} S.~Piccard: 
         Sur un probl\`eme de M.~Ruziewicz de la th\'eorie 
         des relations, 
         {\sl Fundamen\-ta Ma\-thematicae},
         {\bf 29}(1937), 5--9. 

\item{[10]} S. Ruziewicz: 
          Une g\'en\'eralisation d'un th\'eor\`eme de M.~Sierpi\'nski, 
          {\sl Publ. Math.~de l'Universit\'e de Belgrade} 
         {\bf 5} (1936),
         23-27.

\bigskip

\line{\vtop{\hbox{P\'eter Komj\'ath}
            \hbox{Department of Computer Science}
            \hbox{E\"otv\"os University}
            \hbox{Budapest, M\'uzeum krt.~6--8}
            \hbox{1088, Hungary}
            \hbox{e-mail: {\tt kope@cs.elte.hu}}}
            \hfill\vtop{\hbox{S.~Shelah}
                         \hbox{Institute of Mathematics}
                         \hbox{the Hebrew University}
                         \hbox{Jerusalem, Israel}
                         \hbox{e-mail: {\tt shelah@math.huji.ac.il}}}}

\bye